\documentclass[11pt]{amsart}

\title[Mal'cev conditions, graphs 
congruences, tolerances]{Mal'cev conditions from graphs for 
congruences and tolerances}
\author{Paolo Lipparini}
\address{Dipartimento di Matematica, Viale della Ricerc ascientifica,
II Universit\`a di Roma (Tor Vergata),
 ROME 
ITALY}
\thanks{The author has received support from MPI and GNSAGA.
We express our gratitude to G. Cz\'edli for many suggestions.} 

\urladdr{http://www.mat.uniroma2.it/\textasciitilde lipparin}
\newtheorem{Theorem}{Theorem}

\newtheorem{theorem}[Theorem]{Theorem}

\newtheorem{corollary}[Theorem]{Corollary}

\theoremstyle{definition}

\newtheorem{definition}[Theorem]{Definition}

\newtheorem{remarks}[Theorem]{Remarks}

\newcommand{\alg}{\mathbf} 
\def\v{\mathcal V}  

\DeclareMathOperator{\RR}{\mathit{R}}

\begin{document}

\begin{abstract} 
In \cite{acta} we showed that, under certain conditions,
a variety satisfies a given congruence identity if and only if it satisfies the
same tolerance identity. Moreover, we mentioned that a parallel result
holds for Mal'cev conditions arising from graphs. In the present note we
give detailed definitions and state these more general results.
\end{abstract} 

\maketitle

The motivating example for our constructions is the following.
To every $ \{ \circ, \cap\} $-term $p$
one can associate a labeled graph ${\mathbf G}_p $ with 2 distinguished vertices.
See \cite{C1,C2,C3,CD,acta} for details. 
Notice that not every graph can be constructed
as the graph associated with some term, as shown by the example
of the graph with 4 vertices in which every pair of vertices
is connected through some edge.
Our aim here is to show that certain theorems previously stated only
for graphs associated with terms do hold for arbitrary graphs, as
sketched in \cite[Section 7]{acta}.

The classical theory of Mal'cev conditions associates a
 condition $M(p \subseteq  q)$ to each $ \{ \circ, \cap\}$-inclusion
$p \subseteq  q$. However, the definition 
of $M(p \subseteq q)$ is actually given in function 
of the graphs
${\mathbf G}_p$ and
${\mathbf G}_q$ alone,
 with no further reference to the terms
$p$ and $q$. This is particularly evident on \cite[p. 224]{CD}; however, 
it is also implicit in
\cite{C1,C2,C3}.
One can define exactly in the same way
a strong Mal'cev condition
$M({\mathbf G},{\mathbf H})$,
where
${\mathbf G}$
and 
${\mathbf H}$
are arbitrary edge-labeled graphs with
distinguished vertices, no matter
whether or not
${\mathbf G}$ and 
${\mathbf H}$
are associated with terms. 
 We can treat with no additional efforts the slightly more general
case in which we allow more than 2 distinguished vertices.

\begin{definition} \label{MGH} 
Suppose that ${\mathbf G}, {\mathbf H}$
are graphs with edges labeled by 
$n$ labels $ \alpha _1, \dots, \alpha _n$
and
with $h$ distinguished vertices
$d_1,\dots,d_h$ and $e_1,\dots,e_h$.
Let $V$ and
 $W$ denote
 the set of all vertices of  ${\mathbf G}$ and $ {\mathbf H}$.

For each $i$ with $1 \leq i \leq n$,
let $\sim_i$ be the least equivalence relation on $V$ 
such that   $v \sim_i v'$ whenever $v$ and $ v'$
are vertices of $V$ which are connected by some
edge labeled by $ \alpha _i$.
For each $i$ with $1 \leq i \leq n$,
fix $\pi_i$ to be any function from $V$ to an arbitrary set
of variables with the property that $\ker \pi_i=\sim_i$.   

If ${\mathbf G}$ is finite, the 
\emph{Mal'cev condition}
$M({\mathbf G}, {\mathbf H})$ 
involves operations 
$t_w$ ($w \in W$) depending
on $|V|$ variables (in fact, we shall identify 
the variables of $t_w$ with the vertices of 
${\mathbf G}_p$).
Given a fixed arbitrary enumeration
$v_1, \dots, v_m$ of $V$,  
the identities of 
$M({\mathbf G}, {\mathbf H})$ 
are the following:
\[
v_{d_k}=t_{e_k}(v_1,v_2, \dots, v_m)
\] 
for $k=1,\dots,h$, plus
 all the identities:
\[
\tag{m$ _{w, w',i} $ } t_w(\pi_i(v_1), \pi_i(v_2), \dots, \pi_i(v_m))=
t_{w'} (\pi_i(v_1), \pi_i(v_2), \dots, \pi_i(v_m)),
\] 
whenever $w$ and $ w'$ are vertices of  
${\mathbf G}$ connected by an edge labeled
$ \alpha _i$. 
\end{definition} 

In the particular case when 
${\mathbf G}$ and ${\mathbf H}$ are associated with terms, 
we get the classical strong Mal'cev condition associated to an 
inclusion.
Namely, if 
${\mathbf G}={\mathbf G}_p$ and ${\mathbf H}={\mathbf G}_q$,
then 
$M({\mathbf G}_p,{\mathbf G}_q)$ 
turns out to be equal to the classical
$M(p \subseteq q)$.

The Mal'cev conditions introduced in Definition \ref{MGH} are connected
with properties of congruences in algebras by the next definition.

\begin{definition}\label{gr}     
Suppose that ${\mathbf G}$ is a graph with $h>1$ 
distinguished vertices $d_1,\dots,d_h$ and
with edges  labeled by
the set of labels $\{R_1,...,R_n\}$, where the $R_i$'s are symmetric and reflexive binary relations on some set $A$. 

We say     
that the elements
$a_1, \dots, a_h \in A$
can be \emph{connected by the 
graph}
 ${\mathbf G}$
if and only if
there exists some function
$c$ from the set $V$  of vertices 
of  ${\mathbf G}$
to $A$, sending $v\in V$ to $c_v \in A$, 
such that 
$a_1=c _{d_1}, \dots,
a_h=c_{d_h}$, and such that 
whenever the vertices
$v,w \in V$ are connected by an edge labeled by
$R_i$, then $c_v \RR_i c_w$ (cf. 
\cite[Claim 1]{C3},
\cite[Proposition 3.1]{C2} and  
\cite[Proposition 2.1]{acta}).

We define the $h$-ary 
relation 
${\mathbf G} (R_1,...,R_n)$ on $A$ as follows:
$a_1,\dots,a_h \in {\mathbf G} (R_1,...,R_n)$
if and only if  
$a_1, \dots, a_h \in A$
can be connected by the 
graph
 ${\mathbf G}$.

Notice that if 
${\alg A}$ is an algebra and
the $R_i$'s are compatible relations then
${\mathbf G} (R_1,...,R_n)$ is compatible, too, that is,
a subalgebra of $ {\alg A} ^h$.

In the particular case when the graph 
${\mathbf G}= {\mathbf G}_p$
 is associated to the term $p$
we get that 
$(a_1,a_2)\in {\mathbf G}_p(R_1,...,R_n)$
if and only if 
$(a_1,a_2)\in p(R_1 ,..., R_n)$.

Suppose that  ${\mathbf H}$ is another graph of the same type.
Tolerances will be indicated be $ \Theta, \Phi, \dots$,
and congruences by $ \alpha, \beta \dots$.
We say that 
$ {\alg A} $ \emph{satisfies} 
${\mathbf G}\subseteq {\mathbf H}$
\emph{for tolerances},
or simply that
$ {\alg A} $ \emph{satisfies}
${\mathbf G} (\Theta_1, \dots, \Theta_n)\subseteq
 {\mathbf H}(\Theta_1, \dots, \Theta_n)$,
if and only if the inclusion
${\mathbf G} (\Theta_1, \dots, \Theta_n)\subseteq
 {\mathbf H}(\Theta_1, \dots, \Theta_n)$
holds for all $n$-tuples 
$\Theta_1, \dots, \Theta_n$
 of
tolerances on $ {\alg A} $.
We similarly say that 
$ {\alg A} $ \emph{satisfies} 
${\mathbf G}\subseteq {\mathbf H}$
 \emph{for congruences}, or
$ {\alg A} $ \emph{satisfies} 
${\mathbf G} (\alpha _1, \dots, \alpha _n)\subseteq {\mathbf H}(\alpha _1, \dots, \alpha _n )$
when 
the $R_i$'s are only allowed to vary among
 congruences.

Finally, we say that a variety $\v$ 
 \emph{satisfies} 
${\mathbf G}\subseteq {\mathbf H}$
\emph{for tolerances}, or \emph{congruences}
in case every algebra in $\v$
does.
\end{definition}

The classical arguments by Wille and Pixley 
\cite{P,W}
now furnish the following.

\begin{theorem}\label{WP}
For every variety $\v$ and graphs ${\mathbf G}$ and ${\mathbf H}$, 
if ${\mathbf G}$ is finite, then
the following are equivalent.

(i) $\v$ satisfies ${\mathbf G} (\alpha _1, \dots, \alpha _n)\subseteq {\mathbf H}(\alpha _1, \dots, \alpha _n)$.

(ii) The free algebra in $\v$ 
generated by $|V|$ elements satisfies ${\mathbf G} (\alpha _1, \dots, \alpha _n)\subseteq {\mathbf H}(\alpha _1, \dots, \alpha _n)$.

(iii) $\v$ satisfies the Mal'cev condition $M({\mathbf G},{\mathbf H})$.
   \end{theorem}   

\begin{proof}  
(i) $ \Rightarrow  $ (ii) is trivial.
 
(ii) $ \Rightarrow  $ (iii) Consider the free algebra in $\v$
generated by the set $V$ of the vertices of 
${\mathbf G}$ and, for every $i$, let
$ \alpha_i$ be the smallest congruence containing each
equivalence class of $\sim_i$.
The elements
$d_1,\dots,d_h$ can be connected by the 
graph
 ${\mathbf G}$, that is, 
$ d_1,\dots,d_h \in {\mathbf G} (R_1,...,R_n)$,
hence, by assumption,
$ d_1,\dots,d_h \in {\mathbf H} (R_1,...,R_n)$.

If we write down explicitly this last condition,
(letting $e_1=d_1,\dots,e_h=d_h$)
and since we are in a free algebra, we get exactly
terms giving
$M({\mathbf G},{\mathbf H})$.

(iii) $ \Rightarrow  $ (i). Suppose that $ {\alg A} \in \v$
and that $a_1, \dots,a_h$ can be connected by the graph
${\mathbf G}$. This is witnessed by $|V|$ elements
$a_1, \dots,a_m\in A$.
Then the elements
$t_w(a_1, \dots,a_m)$ ($w \in W$)
witness that
$a_1, \dots,a_h$ can be connected by the graph
${\mathbf H}$, by using the
identities in 
$M({\mathbf G},{\mathbf H})$.
\end{proof}

In the particular case when 
${\mathbf G}$ and ${\mathbf H}$ are associated with terms, 
we get back the classical version of the Wille and Pixley Theorem:
for congruences, a variety satisfies $p \subseteq q$ if and only if  it satisfies
$M(p \subseteq q)$, which equals
$M({\mathbf G}_p,{\mathbf G}_q)$. 

We now want to formulate a result which connects 
conditions on congruences to conditions on tolerances, as in \cite{acta}.

In \cite{arx1,acta} we introduced and studied the following notions.

\begin{definition}  \label{repr} 
A tolerance
$\Theta$ of some algebra ${\alg A} $  is {\em representable} 
if and only if there exists a compatible and reflexive
relation $R$ on ${\alg A} $ such that 
$\Theta= R \circ R^-$ (here, $R^-$ denotes the converse of $R$).

A tolerance
$\Theta$ of some algebra ${\alg A} $  is {\em weakly representable} 
if and only if there exists a set $K$ (possibly infinite) and there are
compatible and reflexive
relations $R_k$ ($k \in K$) on ${\alg A} $ such that 
$\Theta= \bigcap _{k \in K} (R_k \circ R_k^-)$. 
\end{definition}

We say that the graph ${\mathbf G}$ 
is \emph{regular} 
if and only if it is finite and
for every $1 \leq i \leq n$ all 
equivalence classes of 
$\sim_i$ have cardinality $\leq 2$.
Thus, if a term $p$ is \emph{regular}
in the sense of \cite{acta}, 
then the labeled graph associated with $p$
is regular.

The methods of \cite{acta} imply the following
statement, which generalizes \cite[Theorem 3.1]{acta}.

\begin{theorem}\label{contolnuo}
For every variety $\v$ and graphs ${\mathbf G}$ and ${\mathbf H}$, 
if ${\mathbf G}$ 
is regular, then the following are equivalent.

(i) $\v$ satisfies ${\mathbf G} (\alpha _1, \dots, \alpha _n)\subseteq {\mathbf H}(\alpha _1, \dots, \alpha _n)$ for congruences.

(ii) $\v$ satisfies ${\mathbf G} (\Theta _1, \dots, \Theta _n)\subseteq {\mathbf H}(\Theta _1, \dots, \Theta _n)$ for representable tolerances.

(iii) $\v$ satisfies ${\mathbf G} (\Theta _1, \dots, \Theta _n)\subseteq {\mathbf H}(\Theta _1, \dots, \Theta _n)$ for weakly representable tolerances.

(iv) $\v$ satisfies ${\mathbf G} (\Theta _1 \circ \Theta _1, \dots, \Theta _n \circ \Theta _n)\subseteq {\mathbf H}(\Theta _1 \circ \Theta _1, \dots, \Theta _n \circ \Theta _n )$ for arbitrary tolerances.
  \end{theorem}   

Theorem \ref{contolnuo}  has the following immediate corollary
(cf. \cite[Corollary 5.1]{acta}).

\begin{corollary}\label{cornuo}
For every variety $\v$ and graphs ${\mathbf G}$ and ${\mathbf H}$, 
if ${\mathbf G}$ 
is regular, then the following are equivalent.

(i) $\v$ satisfies ${\mathbf G} (\alpha _1, \dots, \alpha _n)\subseteq {\mathbf H}(\alpha _1, \dots, \alpha _n)$.

(ii) $\v$ satisfies ${\mathbf G} (\beta_1 \circ \gamma_1 \circ \beta_1, \dots, \beta_n \circ \gamma_n \circ \beta_n)\subseteq {\mathbf H}(\beta_1 \circ \gamma_1 \circ \beta_1, \dots, \beta_n \circ \gamma_n \circ \beta_n)$.

(iii) For every (equivalently, some) odd $m \geq 1$, $\v $ satisfies 
 ${\mathbf G} (\beta_1 \circ_m \gamma_1, \dots, \beta_n \circ_m \gamma_n)\subseteq {\mathbf H}(\beta_1 \circ_m \gamma_1, \dots, \beta_n \circ_m \gamma_n)$.
\end{corollary}

There is a version of Theorem \ref{contolnuo}    
in which the assumption that
${\mathbf G}$ 
is regular is not necessary.
If ${\mathbf G}$ is finite,
compute 
the  integers $k_i$ (which depend only on ${\mathbf G}$) as follows.
If $X$ is an
equivalence class of 
$\sim_i$
 in the graph ${\mathbf G}$
and $x \in X$,
let $k_i(x,X)$
 be the smallest integer such that
every element of $X$ can 
be
connected to $x$ by a path of length $\leq k_i(x,X)$
contained in $X$. 
Let $k_i(X)$ be 
$\inf\{
k_i(x,X)|x \in X
\}$ and finally let $k_i=2\cdot\sup\{
k_i(X)|X \text{\ an\ equivalence\ class\ of\ } 
\sim_i
\}$ (compare the remark after \cite[Proposition 7.6]{acta}).
If $ \Theta $ is a tolerance, let $\Theta^k$
denote $ \Theta \circ \Theta \dots$ with $k$
occurrences of $ \Theta $.

\begin{theorem}\label{contolnuok}
For every variety $\v$ and graphs ${\mathbf G}$ and ${\mathbf H}$, 
the following are equivalent.

(i) $\v$ satisfies ${\mathbf G} (\alpha _1, \dots, \alpha _n)\subseteq {\mathbf H}(\alpha _1, \dots, \alpha _n)$ for congruences.

(ii) $\v$ satisfies ${\mathbf G} (\Theta _1, \dots, \Theta _n)\subseteq {\mathbf H}(\Theta _1^{k_1}, \dots, \Theta _n^{k_n})$ for tolerances. 
  \end{theorem}

\begin{remarks} \label{rmk2}
(a) A generalization of \cite[Remark 7.2]{acta}
along the lines described here holds, too.

(b) We have other applications of the main trick used in the proof
of \cite[Theorem 3.1]{acta}). 
For example, in certain cases, we 
may allow $q$ there to contain the operation
$ \alpha \dot{\cup} \beta =
\overline{ \alpha \cup \beta } $,
where $\overline{X} $ denotes the least
admissible relation generated by $X$.
Moreover, we
have results in which
$p, q$ are allowed to contain
the ternary operations
$K( \alpha , \beta ; \gamma )$ 
and $K _{rs} ( \alpha , \beta ; \gamma )$ 
introduced in 
\cite[Definition 3.4 and p. 867]{lpumi} (cf. also \cite{ctfr}).
 We know still further applications of the arguments in the proof of
\cite[Theorem 3.1]{acta}) and Theorem \ref{contolnuo}.
All such generalizations are not yet published. 
\end{remarks}

\def\cprime{$'$} \def\cprime{$'$}
\providecommand{\bysame}{\leavevmode\hbox to3em{\hrulefill}\thinspace}

\end{document}